\documentclass[oneside,11pt]{amsart}
\usepackage{amsmath, amsfonts,amsthm,times,graphics}
 \makeatletter
\renewcommand*\subjclass[2][2000]{%
  \def\@subjclass{#2}%
  \@ifundefined{subjclassname@#1}{%
    \ClassWarning{\@classname}{Unknown edition (#1) of Mathematics
      Subject Classification; using '1991'.}%
  }{%
    \@xp\let\@xp\subjclassname\csname subjclassname@#1\endcsname
  }%
}
 \makeatother
\newtheorem{theorem}{Theorem}[section]
\newtheorem{lemma}[theorem]{Lemma}
\newtheorem{corollary}[theorem]{Corollary}
\newtheorem{proposition}[theorem]{Proposition}
\theoremstyle{definition}
\newtheorem{definition}[theorem]{Definition}
\newtheorem{example}[theorem]{Example}

\numberwithin{equation}{section}
\newcommand{\abs}[1]{\lvert#1\rvert}

\begin{document}

\title{Lipschitz spaces and harmonic mappings}



\author{David Kalaj}
\address{University of Montenegro, Faculty of Natural Sciences and
Mathematics, Cetinjski put b.b. 81000 Podgorica, Montenegro}
\email{davidk@cg.yu} \subjclass {Primary 30C55, Secondary 31C62}
\begin{abstract}
In \cite{kamz} the author proved that every quasiconformal harmonic
mapping between two Jordan domains with $C^{1,\alpha}$, $0<\alpha\le
1$, boundary is bi-Lipschitz, providing that the domain is convex.
In this paper we avoid the restriction of convexity. More precisely
we prove: any quasiconformal harmonic mapping between two Jordan
domains $\Omega_j$, $j=1,2$, with $C^{j,\alpha}$, $j=1,2$ boundary
is bi-Lipschitz.
\end{abstract} \maketitle

\section{Introduction and notation}
A function $w$ is called \emph{harmonic} in a region $D$ if it has
form  $w=u+iv$ where $u$ and $v$ are real-valued harmonic
functions in $D$. If $D$ is simply-connected, then there are two
analytic functions $g$ and $h$ defined on $D$ such that $w$ has
the representation $$w=g+\overline h.$$

If $w$ is a harmonic univalent function, then by Lewy's theorem
(see \cite{hl}), $w$ has a non-vanishing Jacobian and
consequently, according to the inverse mapping theorem, $w$ is a
diffeomorphism. If $k$ is an analytic function and $w$ is a
harmonic function then $w\circ k$ is harmonic. However $k\circ w$,
in general is not harmonic.

Let $$P(r,x-\varphi)=\frac{1-r^2}{2\pi
(1-2r\cos(x-\varphi)+r^2)}$$ denotes the Poisson kernel. Then
every bounded harmonic function $w$ defined on the unit disc
$\mathbf U:=\{z:|z|<1\}$ has the following representation
\begin{equation}\label{e:POISSON}
w(z)=P[w_b](z)=\int_0^{2\pi}P(r,x-\varphi)w_b(e^{ix})dx,
\end{equation}
where $z=re^{i\varphi}$ and $w_b$ is a bounded integrable function
defined on the unit circle $S^1:=\{z:|z|=1\}$.

Let $A=\begin{pmatrix}
  a_{11} & a_{12} \\
  a_{21} & a_{22}
\end{pmatrix}.$ We will consider two
matrix norms:  $$|A|=\max\{|Az|: z\in \Bbb R^2, |z|=1\}\
\text{and}\ \ |A|_2=({\sum_{i,j}a_{i,j}^2})^{1/2}, $$ and the
matrix function $$l(A)=\min\{|Az|: |z|=1\}.$$ Let $w=u+iv:D\mapsto
G$, $D, G\subset\Bbb C$, be differentiable at $z\in D$. By $\nabla
w(z)$ we denote the matrix $\begin{pmatrix}
  u_{x} & u_{y} \\
  v_{x} & v_{y}
\end{pmatrix}.$ For the matrix $\nabla
w$ we have $$|\nabla w|=|w_z|+|w_{\bar z}|,$$ $$|\nabla
w|_2=(|w_x|^2+|w_y|^2)^{1/2}=\sqrt 2(|w_z|^2+|w_{\bar
z}|^2)^{1/2}$$ and $$l(\nabla w)=||w_z|-|w_{\bar z}||.$$ Thus
\begin{equation}\label{twonorm}|\nabla w|\le |\nabla w|_2\le \sqrt 2 |\nabla
w|.
\end{equation}

A  homeomorphism  $w\colon D\mapsto G, $ where $D$ and $G$ are
subdomains of the complex plane $\mathbf C,$ is said to be
$K$-quasiconformal (K-q.c), $K\ge 1$, if $w$ is absolutely
continuous on a.e. horizontal and a.e. vertical line and
\begin{equation}\label{defqc} \Big|\frac{\partial w}{\partial
x}\Big|^2+ \Big|\frac{\partial w}{\partial y}\Big|^2\le 2
KJ_w\quad \text{a.e. on $D$},\end{equation} where $J_w $ is the
Jacobian of $w$ (cf. \cite{Ahl}, pp. 23--24). Notice that
condition (\ref{defqc}) can be written as $$|w_{\bar z}|\le
k|w_z|\quad \text{a.e. on $D$ where $k=\frac{K-1}{K+1}$ i.e.
$K=\frac{1+k}{1-k}$ },$$ or in its equivalent form
\begin{equation}\label{defqc1} \frac{(|\nabla w|)^2}K \le J_w \le
K (l(\nabla w))^2.
\end{equation}

We will mostly focus on harmonic quasiconformal mappings between
Jordan domains with smooth boundary and will investigate their
Lipschitz character.

Recall that a mapping $w:D\mapsto G$ is said to be $C-$Lipschitz
($C>1$) ($c-$co-Lipschitz) ($0<c$) if $$|w(z_2)-w(z_1)|\le
C|z_2-z_1|, \, z_1,z_2\in D,$$ $$(c|z_2-z_1|\le |w(z_2)-w(z_1)|,
\, z_1, z_2\in D).$$

\section{Background and statement  of the main result} It is well
known that a conformal mapping of the unit disk onto itself has
the form $$w=e^{i\varphi}\frac{z-a}{1-z\bar a},\,
\varphi\in[0,2\pi),\, |a|<1.$$ By the Riemann mapping theorem
there exists a Riemann conformal mapping of the unit disk onto a
Jordan domain $\Omega=\mathrm{int}\ \gamma$. By Caratheodory's
theorem it has a continuous extension to the boundary. Moreover if
$\gamma\in C^{n,\alpha}$, then the Riemann conformal mapping has
$C^{n,\alpha}$ extension to the boundary, see \cite{w2}.
 Conformal mappings are
quasiconformal and harmonic. Hence quasiconformal harmonic mappings
are natural generalization of conformal mappings. The first
characterization of quasiconformal harmonic mappings was started by
O. Martio in \cite{Om}. Hengartner and Schober have shown that, for
a given second dilatation ($a=\overline{f_{\bar z}}/{f_z}$, with
$||a||<1$) there exist a q.c. harmonic mapping $f$ between two
Jordan domains with analytic boundary (\cite[Theorem~4.1]{hs}).
Recently there has been a number of authors who are working on the
topic. Using the result of E. Heinz (\cite{HE}): If $w$ is a
harmonic diffeomorphism of the unit disk onto itself with $w(0)=0$,
then $$|w_z|^2+|w_{\bar z}|^2\ge \frac{1}{\pi^2};$$ O. Martio
(\cite{Om}) observed that, every quasiconformal harmonic mapping of
the unit disk onto itself is co-Lipschitz. Mateljevic, Pavlovic and
Kalaj, have shown that the family of quasiconformal and harmonic
mapping share with conformal mappings the following property: if $w$
is harmonic q.c. mapping of the unit disk onto a Jordan domain with
rectifiable boundary, then $w$ has absolutely continuous extension
to the boundary, see \cite{Kalaj}. What happens if the boundary of a
co-domain is "smoother than rectifiable"?  M. Pavlovic \cite{MP},
proved that every quasiconformal selfmapping of the unit disk is
Lipschitz continuous, using the Mori's theorem on the theory of
quasiconformal mappings. Partyka and Sakan (\cite{pk}) yield
explicit Lipschitz and co-Lipschitz constants depending on constant
of quasiconformality. Since the composition of a harmonic mapping
and of a conformal mapping is itself harmonic, using Kellogg's
theorem (Proposition~\ref{oneone}), these theorems have a
generalization to the class of mappings from arbitrary Jordan domain
with $C^{1,\alpha}$ boundary to the unit disk. However the
composition of a conformal and a harmonic mapping is not, in
general, a harmonic mapping. This means in particular that the
results of this kind for arbitrary co-domain do not follow from the
case of the unit disk and Kellogg's theorem. The situation of
co-domain different from the unit disk firstly has been considered
in \cite{KP}, and there has been shown that every harmonic
quasiconformal mapping of the half-plane onto itself is
bi-Lipschitz. Moreover there have been given two caracterisations of
those mapping, the first one in terms of boundary mapping, using the
Hilbert transformations (\cite{ZY}) and the second one deals with
integral representation, with the help of analytic functions.
Concerning those situations (the disk and the half-plane) see also
\cite{MMM}. The author (\cite{Dk}) extended Heinz theorem
(\cite{HE}) for the harmonic mappings from the unit disk onto a
convex domain. This in turn implies that quasiconformal harmonic
mappings of the unit disk onto a convex domain are co-Lipschitz
(\cite{kalajpub}). Using the new method the results (\cite{MP}) have
been extended properly by the author and Mateljevic in \cite{kamz},
\cite{MM}, and \cite{km}. The extensions are:

Let $\Omega$ and $\Omega_1$ be Jordan domains, let $\mu\in(0,1]$,
and let $f:\Omega_1\mapsto \Omega$ be a harmonic homeomorphism.
Then: (a) If $f$ is $K$ q.c and $\partial \Omega_1, \partial
\Omega\in C^{1,\mu}$, then $f$ is Lipschitz with Lipschitz
constant $c_0(\Omega_1,\Omega,K,w(a))$. Moreover for almost every
$t\in
\partial \Omega_1$ there exists
\begin{equation}\label{nabla}\lim_{z\stackrel {\angle}\to t, z\in
\Omega_1}\nabla{f(z)}=\nabla f(t).\end{equation} (b) If $f$ is q.c
and if $\partial \Omega_1,
\partial \Omega\in C^{1,\mu}$ and $\Omega$ is convex, then $f$
is bi-Lipschitz; (c) If $\Omega_1$ is the unit disk, $\Omega$ is
convex, and $ \partial \Omega\in C^{1,\mu}$, then $f$ is
quasiconformal if and only if its boundary function $f_b$ is
bi-Lipschitz and the Hilbert transformations of its derivative is in
$L^\infty$. (d) If $f$ is q.c and if $\Omega$ is convex then the
boundary functions $f_b$ is bi-Lipschitz in the Euclidean metric and
Cauchy transform $ C[f_b']$  of its derivative is in $L^\infty$. (e)
If $f$ is q.c and if $\Omega$ is convex then the inverse of boundary
functions  $g_b$ is Lipschitz in the Euclidean metric and Cauchy
transform $ C[g_b']$  of its derivative is in $L^\infty$. Concerning
the items (a), (b) and (c) we refer to \cite{kamz}, and for the
items (d) and (e) see \cite{MM} and \cite{napoc}. (f) Let $f$ be a
quasiconformal $C^2$ diffeomorphism from the $C^{1,\alpha}$ Jordan
domain $\Omega_1$ onto the $C^{2,\alpha}$ Jordan domain $\Omega$. If
there exists a constant $M$ such that
\begin{equation}\label{newinequ8}
|\Delta f|\le M|f_z\cdot f_{\bar z}|\,, \quad  z\in
\Omega,\end{equation} then $f$ has bounded partial derivatives. In
particular,  it is a Lipschitz mapping. For the item (f) we refer
to \cite{km}. The result (f) has been generalized in
\cite{kalajmatelj} as follows: (g) Let $f$ be a quasiconformal
$C^2$ diffeomorphism from the plane domain $\Omega_1$ with
$C^{1,\alpha}$ compact boundary onto the plane domain $\Omega$
with $C^{2,\alpha}$ compact boundary. If there exist constants $M$
and $N$ such that
\begin{equation}\label{newinequ2}
|\Delta f|\le M|\nabla f|^2+N\,, \quad  z\in \Omega
,\end{equation} then $f$ has bounded partial derivatives in
$\Omega_1$. In particular it is a Lipschitz mapping in $\Omega_1$.

 For several dimensional generalizations we refer to \cite{mz},
\cite{mv} and \cite{akm}.

Because of the lack of generalization of the Heinz theorem for non
convex domains, it was intrigue to investigate the q.c. harmonic
mappings of the unit disk onto the image domain that is not
convex. Namely it has been an open problem until now that, if the
assumption of convexity on an image domain $\Omega$ was important
or not in proving the theorem that a harmonic q.c. mapping of the
unit disk onto $\Omega$ is bi-Lipschitz.

In the following theorem we avoid the restriction of convexity.

\begin{theorem}[The main theorem]\label{krye}
Let $w=f(z)$ be a $K$ quasiconformal harmonic mapping between a
 Jordan domain $\Omega_1$ with $C^{1,\alpha}$ boundary and a
Jordan domain $\Omega$ with $C^{2,\alpha}$ boundary. Let in addition
$b\in \Omega_1$ and $a=f(b)$. Then $w$ is bi-Lipschitz. Moreover
there exists a positive constant $c=c(K,\Omega,\Omega_1, a, b)\ge 1$
such that
\begin{equation}\label{pocetna} \frac 1c |z_1-z_2|\le
|f(z_1)-f(z_2)|\le c|z_1-z_2|,\,\,\,\,z_1,z_2\in
\Omega_1.\end{equation}
\end{theorem}

\section{The proof of the main theorem}

The key of the proof is Lemma~\ref{hopf}, which could be
considered as a global version of the following well known lemma:

\begin{lemma}[Hopf's Boundary Point lemma]\cite{prot} and \cite{hopf}. Let $u$ satisfies $\Delta u\ge 0$ in $D$ and $u\le M$
in D, $u (P) = M$ for some $P\in \partial D$. Assume that $P$ lies
on the boundary of a ball $B\subset D$. If $u$ is continuous on
$D\cup P$ and if the outward directional derivative
$\frac{\partial u}{\partial n}$ exists at $P$, then $u\equiv M$
 or $$\frac{\partial u}{\partial n}>0.$$
\end{lemma}

\begin{lemma}\label{hopf} Let $u$ satisfies $\Delta u\ge 0$ in
$R_\varrho=\{z:\varrho\le |z|<1\}$, $0<\varrho<1$, $u$ be
continuous on $\overline{R_\varrho}$, $u< 0$ in $R_\varrho$, $u
(t) = 0$ for $t\in S^1$. Assume that the radial derivative
$\frac{\partial u}{\partial r}$ exists almost everywhere at $t\in
S^1$. Let $M(u,\varrho):=\max_{|z|=\varrho}u(z)$. Then for the
positive constant
\begin{equation}\label{ndih}c(u,\varrho)=\frac{2M(u,\varrho)}{\varrho^2(1-e^{1/\varrho^2-1})}\end{equation}
there holds \begin{equation}\label{ndihma}\frac{\partial
u(t)}{\partial r}>c(u,\varrho), \text{ for a.e. } t\in
S^1.\end{equation}
\end{lemma}

\begin{proof}

Consider the auxiliary function $h^A_\varrho(z) = e^{-A|z|^2}-e^{
-A}$, where $A>0$ is a constant to be chosen later. Then $$ \Delta
h^A_\varrho(z) = 4Ae^{A|z|^2}(A|z|^2-1). $$

Hence it has the property that $h^A_\varrho(z)>0$, $z\in
R_\varrho$, and that
\begin{equation}\label{2}\Delta h^A_\varrho
\ge 0, \varrho \le |z|\le 1,\end{equation} if
\begin{equation}\label{cona}A\ge {\varrho}^{-2}, \text{for example $A=\varrho^{-2}$}.\end{equation} The function
$h^A_\varrho(z)$ is of class $C^2 $ in $R_\varrho$, and
\begin{equation}\label{3} h^A_\varrho(z)=0\ \  \text{on} \ \
S^1.\end{equation}

The function $v^A_\varrho = u+\varepsilon h^A_\varrho(z)$,
$\varepsilon
> 0,$ is of class $C^2$ in the interior of $R_\varrho$ and continuous in $R_\varrho$.
Moreover, by \eqref{3}, \begin{equation}\label{4} v^A_\varrho \le
0\text{ on }S^1.\end{equation}

As $M(u,\varrho)<0$ we can choose a constant $\varepsilon$ so that
$$M(u,\varrho)+\varepsilon (e^{-A\varrho^2}-e^{ -A})\le 0.$$ For
example
\begin{equation}\label{eps}\varepsilon =
\frac{M(u,\varrho)}{e^{-A}-e^{-A\varrho^2}}.\end{equation} Then we
have
\begin{equation}\label{5} v^A_\varrho \le 0 \ \ \text{also on}\ \ S(0,\varrho)
.\end{equation}

By the hypothesis, $\Delta u \ge 0$ in $R_\varrho$, and by
\eqref{2} it follows
\begin{equation}\label{6} \Delta v_\varrho^A > 0, z\in R_\varrho.\end{equation} \eqref{4}, \eqref{5}, and \eqref{6}
imply that $v_\varrho^A\le 0$ holds in the whole of $R_\varrho$.
This follows from the elementary fact that $v_\varrho^A$ cannot
have a positive maximum in the interior of $R_\varrho$. But
$v_\varrho^A \le 0$ in $R_\varrho$ and $v_\varrho^A = 0$ at $t\in
S^1$ implies that

$$0\le \lim_{R\to 1-0}\frac{v_\varrho^A(Rt)-v_\varrho^A(t)}{R-1}=
\frac{\partial v_\varrho^A(t)}{\partial r}=\frac{\partial
u(t)}{\partial r}+\varepsilon\frac{\partial
h_\varrho^A(t)}{\partial r}.$$ Furthermore
$$\displaystyle\min_{s\in S^1}\frac{\partial
h_\varrho^A(s)}{\partial r}=-2Ae^{-A}<0.$$ Thus for almost every
$t\in S^1$ there holds
\begin{equation}\frac{\partial u(t)}{\partial
r}\ge-\varepsilon\min_{s\in S^1}\frac{\partial
h_\varrho^A(s)}{\partial
r}=\frac{2AM(u,\varrho)}{1-e^{(1-\varrho^2)A}}=:c(u,\varrho)>0.\end{equation}
\end{proof}

To continue we need the following propositions:

\begin{proposition}[Kellogg]\cite{G}\label{oneone} {\it If a domain
$D=Int(\Gamma)$ is $ C^{1,\alpha}$ and $\omega$ is a conformal
mapping of $\mathbb{U}$ onto $D$, then $\omega'$ and $\ln \omega'$
are in $Lip_\alpha$. In particular, $|\omega'|$ is bounded from
above and below on ${\mathbf U}$} by two positive constants.
\end{proposition}

Let $\Gamma$ be a smooth Jordan curve and $\beta (s)$ the angle of
the tangent as a function of arc length. We say that $\Gamma$ has
a Dini-continuous curvature if $\beta'(s)$ is continuous and $$
|\beta'(s_2)-\beta'(s_1)|\le \omega_1(s_2-s_1)\,\, (s_1<s_2),$$
where $\omega_1(x)$ is an increasing function that satisfies
$$\int_0^1\frac{\omega_1(s)}{s}ds<\infty.$$ The next proposition
is due to Kellogg and to Warschawski.
\begin{proposition} \cite[Theorem~3.6]{w}\label{onetwo}. Let $\omega$ be a conformal mapping of the unit disk onto a Jordan domain that is bounded by a Jordan curve with Dini- continuous
curvature. Then $\omega'' (z)$ has a continuous extension to
$\overline{\Bbb U}$. In particular $|\omega''|$ is bounded from
above on $\Bbb U$.
\end{proposition}

Notice that if $\Gamma$ is $C^{2,\alpha}$ then $\Gamma$ has
Dini-continuous curvature. We will finish the proof of
Theorem~\ref{krye} using the following lemma.

\begin{lemma}\label{lemica}
Let $w=f(z)$ be a $K$ quasiconformal harmonic mapping of the unit
disk onto a $C^{2,\alpha}$ Jordan domain $\Omega$ such that
$w(0)=a\in \Omega $. Then there exists a constant
$C(K,\Omega,a)>0$ such that $$\abs{\frac{\partial w}{\partial
r}(t)}\ge C(K,\Omega,a)\text{ for almost every } t\in S^1.$$
\end{lemma}

\begin{proof}
Let $g$ be a conformal mapping of $\Omega$ onto the unit disk with
$g(a)=0$. Take $w_1=g\circ w$. Then
\begin{equation}\label{ew}\begin{split}\Delta w_1 &= 4g''(w) w_z\cdot
w_{\bar z} +g'(w)\Delta w\\&=4g''(w) w_z\cdot w_{\bar
z}=4\frac{g''}{|g'|^2}{w_1}_z\cdot {w_1}_{\bar
z}.\end{split}\end{equation}

Combining \eqref{ew} and \eqref{defqc1} we obtain
\begin{equation}\label{kompo} |\Delta w_1|\le
\frac{|g''|}{|g'|^2}|(|\nabla w_1|^2-l(\nabla w_1)^2) \le
\left(1-\frac{1}{K^2}\right)\frac{|g''|}{|g'|^2}|\nabla
w_1|^2.\end{equation} Let $h(z)=|w_1|^2$. Let us find two
constants $B>0$ and $\varrho\in(0,1)$ such that the function
$$\varphi(z):=\chi(h(z))= \frac{1}{B}(e^{Bh(z)}-e^{B})$$
 is subharmonic on
$\{z: \varrho<|z|<1\}$. Clearly $\varphi(z)\le 0$. On the other
hand we have

\begin{equation}\label{kompi}\Delta \varphi =\chi''(h)|\nabla h|^2+\chi'(h)\Delta
h.\end{equation} Furthermore
\begin{equation}\label{delc}\Delta h=2|\nabla w_1|_2^2+
2\left<\Delta w_1,w_1\right> .\end{equation} Let $w_1=\rho s$,
$\rho=|w_1|$, $s=e^{i\psi}$. Then

\begin{equation}\label{delca}|\nabla h|=2\rho|\nabla
\rho|.\end{equation} To continue observe that $$\nabla
w_1=(\nabla\rho)^t s+\rho \nabla s$$ and thus $$|\nabla w_1\,
l|^2=|\rho \nabla s\, l|^2+|\nabla \rho\, l \cdot
s|^2+2\rho\nabla\rho l\left<\nabla s\,l,s\right>, \, l\in \Bbb
R^2.$$ Hence
\begin{equation}\label{helpeq}|\nabla w_1\,l|^2=\rho^2|\nabla s\,l|^2+|\nabla \rho
\,l|^2.\end{equation} Choose $l_1: |l_1|=1$ so that $\nabla s
l_1=0$. Then  by (\ref{helpeq}) we infer $$|\nabla
w_1\,l_1|\le|\nabla \rho \,l_1|.$$ According to the definition of
quasiconformal mappings we obtain
\begin{equation}\label{quas}K^{-1}|\nabla w_1|\le |\nabla
\rho|.\end{equation}

From \eqref{delca} and \eqref{quas} it follows that

\begin{equation}\label{qu}
|\nabla h|\ge \frac{2\rho}{K}|\nabla w_1|.
\end{equation}

Combining \eqref{twonorm}, \eqref{kompo},
\eqref{kompi}, \eqref{delc} and \eqref{qu} we obtain

\begin{equation}\label{gati}
\Delta \varphi\ge
\left(\chi''\frac{4\rho^2}{K^2}+2\chi'-2\left(1-\frac{1}{K^2}\right)\chi'\frac{|g''|}{|g'|^2}\right)|\nabla
w_1|^2.
\end{equation}

Furthermore
\begin{equation}\label{firstder} \chi'(h)=e^{B h}
\end{equation}
and
\begin{equation}\label{secder} \chi''(h)=Be^{Bh}.
\end{equation}
By \eqref{gati}, \eqref{firstder} and \eqref{secder} we obtain

\begin{equation}\label{A}\Delta \varphi\ge \left(B\frac{4\rho^2}{K^2}+2-2\left(1-\frac{1}{K^2}\right)
\frac{|g''|}{|g'|^2}\right)e^{Bh(z)}|\nabla
w_1|^2.\end{equation}

As $w_1=\rho s$ is $K$ quasiconformal selfmapping of the unit disk
with $w_1(0)=0$, by Mori's theorem (\cite{wang}) it satisfies the
doubly inequality:
\begin{equation}\label{mori}
\left|\frac z{4^{1-1/K}}\right|^K \le \rho\le 4^{1-1/K} |z|^{1/K}.
\end{equation}

By \eqref{mori} for $\varrho\le |z|\le 1$ where
\begin{equation}\label{varh}\varrho:={4^{-K}}\end{equation} we have
\begin{equation}\label{rhoeq}
\rho\ge 4^{1-K^2-K}.
\end{equation}
Now we choose $B$ such that
$$\frac{4B\rho^2}{K^2}+2-2\left(1-\frac{1}{K^2}\right)\frac{|g''|}{|g'|^2}\ge
0,$$ i.e. in view of Propositions~\ref{oneone} and \ref{onetwo},
and \eqref{rhoeq}, for example take:

\begin{equation}\label{eqb}B:=\max\left\{\frac{1}{2}\sup_{z\in
\Omega}\left|1-\left(1-\frac{1}{K^2}\right)\frac{|g''|}{|g'|^2}\right|K^24^{K^2+K-1},1\right\}.\end{equation}

According to Lemma~\ref{hopf}, and to \eqref{nabla} the function
$$\varphi(z)=\chi(h(z))= \frac{1}{B}(e^{Bh(z)}-e^B)$$ satisfies

$$\frac{\partial \varphi}{\partial
R}(t)=e^{Bh(t)}\left<g'(w(t))\cdot \frac{\partial w}{\partial
R}(t), w_1(t)\right>\ge c(\varphi,\varrho),$$ almost everywhere in
$S^1$, where $c(\varphi,\varrho)$ is defined by \eqref{ndih}. On
the other hand by the right hand inequality in \eqref{mori} it
follows that
\begin{equation}\label{ocenaf}
\varphi(z)\le \frac{1}{B}(e^{4^{-\frac 2K}B}-e^B) \text{ for $|z|
= \varrho$.}
\end{equation}
Thus
\begin{equation}\label{minvar}M(\varphi,\varrho)=\max_{|z|=\varrho}
\varphi(z) \le \frac{1}{B}(e^{4^{-\frac
2K}B}-e^B)<0.\end{equation} According to \eqref{ndih} and
\eqref{ndihma} it follows that $$\abs{\frac{\partial w}{\partial
r}(t)}\ge \frac{e^{-B}c(\varphi,\varrho)}{\max\{|g'(\zeta)|:
\zeta\in\partial\Omega\}}=\frac{2e^{-B}M(\varphi,\varrho)}{\varrho^2(1-e^{1/\varrho^2-1})|g'|_\infty}>0,$$
almost everywhere in $S^1$. By \eqref{varh}, \eqref{eqb} and
\eqref{minvar}, we can take
$$C(K,\Omega,a)=\frac{2e^{-B}M(\varphi,\varrho)}{\varrho^2(1-e^{1/\varrho^2-1})|g'|_\infty}
$$ ($C(K,\Omega,a)$ do not depends on $w=f(z)$).

\end{proof}
\begin{proof}[\bf Proof of Theorem~\ref{krye}]
In view of item a) from the Background of this paper, it is enough
to prove that, $w$ is co-Lipschitz continuous (under the above
conditions). Moreover by Proposition~\ref{oneone} the unit disk
could be taken as the domain of the mapping.

We will consider two cases:
\\
{\bf 1.  CASE "$w\in C^{1}(\overline {\mathbf U})$".}

Let $l(\nabla w)(t)=||w_z(t)|-|w_{\bar z}(t)||$. As $w$ is $K$
q.c., according to Lemma~\ref{lemica} we have
\begin{equation}\label{due}l(\nabla w)(t)\ge\frac{|\nabla w(t)|}{K}\ge
\frac{\abs{\frac{\partial w}{\partial r}(t)}}{K}\ge
\frac{C(K,\Omega,a)}{K},\end{equation} for $t\in S^1$.

Since $w$ is a harmonic diffeomorphism, by the Lewy theorem
(\cite{hl}) ($|w_z|> 0$), it defines the bounded subharmonic
function
\begin{equation}\label{ss}S(z):=\left|\frac{\overline{w_{\bar z}}}{w_z}\right| +
\left|\frac{1}{w_z}\cdot
\frac{C(K,\Omega,a)}{K}\right|\end{equation} on the unit disk.
According to \eqref{due}, $S(z)$ is bounded on the unit circle by
1. By the maximum principle, this implies that $S$ is  bounded on the whole unit disk by 1.

This in turn implies that for every $z\in {\mathbf U}$
\begin{equation}\label{maj}l(\nabla w)(z) \ge  \frac{C(K,\Omega,a)}{K} .\end{equation}
\\
{\bf 2. CASE "$w
 \notin C^{1}(\overline {\mathbf U})$".}

\begin{definition}
{\it Let $G$ be a domain in $\Bbb C$ and let $a\in\partial G$. We
will say that $G_a\subset G$ is a neighborhood of $a$ if there
exists a disk $D(a,r):=\{z:|z-a|<r\}$ such that $D(a,r)\cap
G\subset G_a$.}
\end{definition}
 Let $t=e^{i\beta}\in S^1$, then $w(t)\in
 \partial \Omega$. Let $\gamma$ be an arch length parametrization of $\partial\Omega$
with $\gamma(s)=w(t)$. Since $\partial \Omega \in C^{2,\alpha}$
there exists a neighborhood $\Omega_t$ of $w(t)$ with
$C^{2,\alpha}$ Jordan boundary such that,
\begin{equation}\label{subdomains}\Omega^\tau _t:=\Omega_t +
i\gamma'(s)\cdot \tau\subset \Omega,\text{ and } \partial
\Omega^\tau_t\subset \Omega \text{ for $0< \tau\le \tau_t$\,
($\tau_t>0$) }.\end{equation} An example of a family
$\Omega^\tau_t$ such that $\partial\Omega^\tau_t \in C^{1,\alpha}$
and with the property \eqref{subdomains} has been given in
\cite{kamz}. An easily modification yields a family of Jordan
domains $\Omega^\tau_t$ with $\partial\Omega^\tau_t\in
C^{2,\alpha}$, $0\le \tau\le \tau_t$ with the property
\eqref{subdomains}.

Let $a_t\in \Omega_t$ be arbitrary. Then $a_t+i\gamma'(s)\cdot
\tau\in \Omega_t^\tau$. Take $U_\tau = f^{-1}(\Omega_t^\tau)$. Let
$\eta_t^\tau $ be a conformal mapping of the unit disk onto $U_\tau$
such that $\eta_t^\tau(0)=f^{-1}(a_t +i\gamma'(s)\cdot \tau)$, and
$\arg\frac{ d  \eta_t^\tau}{ d  z}(0)=0$. Then the mapping
$$f_t^\tau(z) := f(\eta_t^\tau(z))-i\gamma'(s)\cdot \tau$$ is a
harmonic $K$ quasiconformal mapping of the unit disk onto $\Omega_t$
satisfying the condition $f_t^\tau(0)= a_t$. Moreover
$$f_t^\tau\in C^1(\overline {\mathbf U}).$$ Using the {\bf CASE
"$w\in C^1(\overline {\mathbf U})$"} it follows that $$|\nabla
f_t^\tau (z)|\ge C(K,\Omega_t, a_t).$$ On the other hand
$$\lim_{\tau\to 0+}\nabla f_t^\tau(z) = \nabla (f\circ \eta_t)(z)$$
on the compact sets of ${\mathbf U}$ as well as $$\lim_{\tau \to
0+}\frac{ d  \eta_t^\tau}{ d  z}(z)= \frac{ d \eta_t}{ d z}(z),$$
where $\eta_t$ is a conformal mapping of the unit disk onto
$U_0=f^{-1}(\Omega_t)$ with $\eta_t(0)=f^{-1}(a_t)$. It follows that
$$|\nabla f_t(z)|\ge C(K,\Omega_t, a_t).$$

Using the Schwartz's reflexion principle to the mapping $\eta_t$,
and using the formula $$\nabla (f\circ \eta_t)(z) = \nabla f\cdot
\frac{ d  \eta_t}{ d  z}(z)$$ it follows that in some neighborhood
$\tilde U_t$ of $t\in S^1$ ($D(t,r_t)\cap \mathbf U\subset \tilde
U_t$ for some $r_t>0$) the function $f$ satisfies the inequality
\begin{equation}\label{local} |\nabla f(z)|\ge \frac{C(K,\Omega_t,
a_t)}{\min\{|\eta_t(\zeta)| : \zeta\in
\partial \tilde U_t\cap S^1\}}=:\tilde C(K,\Omega_t,a_t)>0.
\end{equation}

Since $S^1$ is a compact set it can be covered by a finite family
$\partial\tilde U_{t_j}\cap S^1\cap D(t,r_t/2)$, $j=1,\dots,m$. It
follows that the inequality
\begin{equation}\label{local} |\nabla f(z)|\ge
\min\{\tilde C(K,\Omega_{t_j}, a_{t_j}): j=1,\dots, m\}=:\tilde
C(K,\Omega,a)>0,
\end{equation} there holds in the annulus $$\tilde R=\left\{z:1-\frac{\sqrt 3}{2}\min_{1\le j\le m}r_{t_j}<|z|<1\right\}
\subset \bigcup_{j=1}^m \tilde U_{t_j}.$$
 This implies that the subharmonic function $S$ defined in
\eqref{ss} is bounded in $\mathbf U$. According to the maximum
principle it is bounded by $1$ in the whole unit disk. This in
turn implies again \eqref{maj} and consequently

$$\frac{C(K,\Omega,a)}{K}|z_1-z_2|\le |w(z_1)-w(z_2)|,\ \ \ \
z_1, z_2\in \mathbf U.$$

\end{proof}

\begin{corollary}\label{coro}
If $w$ is q.c. harmonic mapping of the unit disk onto a
$C^{2,\alpha}$ Jordan domain $\Omega$, then
$\mathrm{ess\,sup}\,\{J_w(z), z\in \Bbb U\}>0$. Recall that by
$J_w$ we denote the Jacobian of $w$.
\end{corollary}

\begin{example}
$w=P[e^{i(x+\sin x)}](z)$, $z\in \Bbb U$ is a harmonic
diffeomorphism of the unit disk onto itself having smooth
extension to the boundary and \[\begin{split}0&\le J_w(-1)\le
|\frac{\partial w}{\partial r}(re^{i\varphi})|_{r=1,
\varphi=\pi}|\cdot|\frac{\partial
w}{\partial\varphi}(e^{i\varphi})|_{
\varphi=\pi}\\&=|\frac{\partial w}{\partial
r}(re^{i\varphi})|_{r=1,
\varphi=\pi}\cdot|(1+\cos\varphi|_{\varphi=\pi})|=0,\end{split}\]
i.e. $J_w(-1)=0$ . Hence the condition of quasiconformality in
Corollary~\ref{coro} is essential.
\end{example}

\subsection{Remarks}
It seems natural that the assumption $\partial \Omega\in
C^{2,\alpha}$ in the main theorem can be replaced by $\partial
\Omega\in C^{1,\alpha}$ however we do not have the proof of this
fact. It remains an open problem, whether the norm of the first
derivative of harmonic diffeomorphism between the unit disk and a
smooth Jordan domain $\Omega$ is bounded bellow by a constant
depending on $\Omega$. The result of this kind was proved by E.
Heinz, \cite{HE},  for the case of $\Omega$ being the unit disk
and by the author in \cite{Dk} for $\Omega$ being a convex domain.
In this paper it was proved that the result hold for harmonic
quasiconformal mappings without the restriction on convexity of
co-domain.

\subsection*{Acknowledgment} I thank the referee for useful comments and suggestions
related to this paper.


\begin{thebibliography}{1}
\bibitem{Ahl} L. Ahlfors: \textit{Lectures on Quasiconformal mappings,}
Van Nostrand Mathematical Studies, D. Van Nostrand 1966.
\bibitem{akm}
M. Arsenovic, V. Kojic and M. Mateljevic: {\it On lipschitz
continuity of harmonic quasiregular maps on the unit ball in $\Bbb
R^n$}, Ann. Acad. Sci. Fenn., Math. Vol 33, 315-318, (2008).
\bibitem{G}
G. M. Goluzin: {\it Geometric function theory }, Nauka Moskva 1966
(Russian).
\bibitem{hs}
W. Hengartner, G. Schober: {\it Harmonic mappings with given
dilatation.}  J. London Math. Soc. (2)  33  (1986),  no. 3,
473--483.
\bibitem{HE}
E. Heinz: {\it On one-to-one harmonic mappings. } Pac. J. Math. 9,
101-105 (1959).
\bibitem{hopf} E. Hopf: {\it A remark on linear elliptic
differential equations of second order}, Proc. Amer. Math. Soc.,
3, 791-793  (1952).
\bibitem{kalajpub}
D. Kalaj: \emph{Quasiconformal harmonic functions between convex
domains}, Publ. Inst. Math., Nouv. Ser. 76(90), 3-20 (2004).
\bibitem{km} D. Kalaj, M. Mateljevi\'c:
{\it Inner estimate and quasiconformal harmonic maps between
smooth domains}, Journal d'Analise Math. 100. 117-132, (2006).
\bibitem{KP}D. Kalaj and M. Pavlovi\'c:
{\it Boundary correspondence under harmonic quasiconformal
homeomorfisms of a half-plane,} Ann. Acad. Sci. Fenn., Math. 30,
No.1, (2005) 159-165.
\bibitem{mz} D. Kalaj: {\it On harmonic quasiconformal self-mappings of the unit
ball}, Ann. Acad. Sci. Fenn., Math. Vol 33, 1-11, (2008).
\bibitem{kamz} D. Kalaj: {\it Quasiconformal harmonic mapping
between Jordan domains}  Math. Z. Volume 260, Number 2, 237-252,
2008.
\bibitem{kalajmatelj}
D. Kalaj, M. Mateljevic: {\it On certain nonlinear elliptic PDE
and quasiconfomal mapps between Euclidean surfaces},
arXiv:0804.2785.
\bibitem{Dk}
D. Kalaj: {\it On harmonic diffeomorphisms of the unit disc onto a
convex domain.} Complex Variables, Theory Appl. 48, No.2, 175-187
(2003).
\bibitem{Kalaj} D. Kalaj: \textit{Harmonic functions and harmonic quasiconformal mappings
between convex domains,} Thesis, Beograd 2002.
\bibitem{MMM}
M. Knezevic, M. Mateljevic: {\it On the quasi-isometries of
harmonic quasiconformal mappings} Journal of Mathematical Analysis
and Applications, 2007; 334 (1) 404-413.
\bibitem{hl}
H. Lewy: {\it On the non-vanishing of the Jacobian in certain in
one-to-one mappings,} Bull. Amer. Math. Soc. 42. (1936), 689-692.
\bibitem{Om}
O. Martio: {\it On harmonic quasiconformal mappings}, Ann. Acad.
Sci. Fenn., Ser. A I 425 (1968), 3-10.
\bibitem{MM} M. Mateljevic:
{\it On quasiconformal harmonic mappings}, unpublished manuscript,
2006.
\bibitem{mv}
M. Mateljevic, M. Vuorinen: {\it On harmonic quasiconformal
quasi-isometries}, arXiv: 0709.4546v1.
\bibitem{napoc} M.  Mateljevi\'c,  {\it Distortion  of harmonic
functions and harmonic quasiconformal  quasi-isometry},   Revue
Roum. Math. Pures Appl.  Vol.{\bf 51},(2006)), 5-6, 711-722
\bibitem{pk} D. Partyka and K. Sakan:
{\it On bi-Lipschitz type inequalities for quasiconformal harmonic
mappings,} Ann. Acad. Sci. Fenn. Math.. Vol 32, pp. 579-594
(2007).
\bibitem{MP} M. Pavlovi\' c:
{\it Boundary correspondence under harmonic quasiconformal
homeomorfisms of the unit disc}, Ann. Acad. Sci. Fenn., Vol 27,
(2002) 365-372.
\bibitem{w} C. Pommerenke: {\it Boundary behavour of conformal maps,}
Springer-Verlag, New York, 1991.
\bibitem{pw}
C. Pommerenke and S.E. Warschawski: {\it  On the quantitative
boundary behavior of conformal maps.} Comment. Math. Helv. 57,
107-129 (1982).
\bibitem{prot} M. H. Protter and H. F. Weinberger: {\it Maximum principles in
differential equations,} Prentice Hall, Englewood Cliks N.J.,
1967.
\bibitem{wang}
C. Wang: {A sharp form of Mori's theorem on Q-mappings,} Kexue
Jilu, 4 (1960), 334-337.
\bibitem{w1} S. E. Warschawski: {\it On differentiability at the boundary in
conformal mapping,} Proc. Amer. Math. Soc, 12 (1961), 614-620.
\bibitem{w2} S. E. Warschawski: {\it On the higher derivatives at the boundary in conformal
mapping,} Trans. Amer. Math. Soc, 38, No. 2 (1935), 310-340.
\bibitem{ZY}
A. Zygmund:  {\it Trigonometric Series I.} Cambrige University
Press, 1958.
\end{thebibliography}
\end{document}